\documentclass[11pt,leqno]{amsart}
\usepackage{amssymb,verbatim,enumerate,ifthen,hyperref}
\usepackage{mathtools} 

\usepackage{cite}
\usepackage{tikz-network}
\usepackage{graphicx}
\usepackage[mathscr]{eucal}
\usepackage[utf8]{inputenc}
\usepackage[T1]{fontenc}
\def\N{\mathbb{N}}
\def\R{\mathbb{R}}

\def\A{\mathscr{A}}

\long\def\comment#1{}

\newtheorem{theorem}{Theorem}[section]
\newtheorem*{theorem*}{Theorem}
\def\Thm#1#2{\ifthenelse{\equal{#1}{*}}{\begin{theorem*}#2\end{theorem*}}
             {\begin{theorem}\label{T#1}#2\end{theorem}}}
\newtheorem{Atheorem}{Theorem}

\def\thm#1{Theorem~\ref{T#1}}
\newtheorem{proposition}[theorem]{Proposition}
\newtheorem*{proposition*}{Proposition}
\def\Prp#1#2{\ifthenelse{\equal{#1}{*}}{\begin{proposition*}#2\end{proposition*}}
{\begin{proposition}\label{P#1}#2\end{proposition}}}
\def\prp#1{Proposition~\ref{P#1}}

\newtheorem{corollary}[theorem]{Corollary}
\newtheorem*{corollary*}{Corollary}
\def\Cor#1#2{\ifthenelse{\equal{#1}{*}}{\begin{corollary*}#2\end{corollary*}}
             {\begin{corollary}\label{C#1}#2\end{corollary}}}
\def\cor#1{Corollary~\ref{C#1}}

\newtheorem{lemma}[theorem]{Lemma}
\newtheorem*{lemma*}{Lemma}
\def\Lem#1#2{\ifthenelse{\equal{#1}{*}}{\begin{lemma*}#2\end{lemma*}}
             {\begin{lemma}\label{L#1}#2\end{lemma}}}

\theoremstyle{definition}
\newtheorem{remark}[theorem]{Remark}
\newtheorem*{remark*}{Remark}
\def\Rem#1#2{\ifthenelse{\equal{#1}{*}}{\begin{remark}\rm #2\end{remark}}
             {\begin{remark}\label{R#1}\rm #2\end{remark}}}

\newtheorem{example}[theorem]{Example}
\newtheorem*{example*}{Example}
\def\Exa#1#2{\ifthenelse{\equal{#1}{*}}{\begin{example*}\rm #2\end{example*}}
             {\begin{example}\label{Ex#1}\rm #2\end{example}}}

\def\eq#1{{\rm(\ref{E#1})}}
\def\Eq#1#2{\ifthenelse{\equal{#1}{*}}
  {\begin{equation*}\begin{aligned}#2\end{aligned}\end{equation*}}
  {\begin{equation}\begin{aligned}\label{E#1}#2\end{aligned}\end{equation}}}

\begin{document}
\begin{flushright}
\end{flushright}
\vspace{5mm}

\date{\today}

\title[On Weighted Convex Graphs]
{On Weighted Convex Graphs}

\author[A. R. Goswami]{Angshuman R. Goswami}
\address[A. R. Goswami]{Department of Mathematics, University of Pannonia,
H-8200 Veszpr\'em, Hungary}
\email{goswami.angshuman.robin@mik.uni-pannon.hu}

\subjclass[2000]{Primary: 05C05, 05C22, 39B82; Secondary: 05C90, 39A12, 52A01}
\keywords{Weighted convexity on graphs;  Convex Sequence; Sequence Embedding, Ulam Stability}

\thanks{The research of the first author was supported by the 
}

\begin{abstract}
The main objective of this paper is to develop Krein-Milman-type theorems and 
Ulam-type stability results for graphs. To establish these results, we introduce several meaningful definitions of vertex-weighted convex graphs inspired by the concept of sequential convexity. We also present a close relationship between the two discrete structures, namely sequential convexity and perfect binary trees. We show that if a graph satisfies a certain convexity property approximately, then this property can be made exact by minimally perturbing the weights assigned to its vertices. Furthermore, we study several structural characterisations, formulate convex minorants for weighted graphs, and derive sandwich-type results. Special emphasis is placed on trees, and an investigation of extremal value problems is also carried out.\\

Various definitions, research backgrounds, motivations, 
and other crucial details are discussed in the following section.
\end{abstract}

\maketitle
\section*{Introduction}
In classical graph theory, the notion of convexity is defined as follows:
Let $G(V,E)$ be a simple connected graph. A subset $U \subseteq V$ is called \emph{convex} if for every pair of vertices $u_1,u_2 \in U$, every vertex that lies on a shortest $u_1$--$u_2$ path in $G$ also belongs to $U$. 
Equivalently, $U$ is convex if
$$
   \mbox{for all}\quad u_1,u_2 \in U, \quad I(u_1,u_2) \subseteq U,
$$
where $I(u_1,u_2)$ denotes the set of all vertices that lie on some shortest path between $u_1$ and $u_2$ in $G$. The classical notion of a convex set inspires this definition. \\

Researchers have defined various similar notions for convexity and have extensively studied several combinatorial properties of graphs. Carathéodory, Helly, and 
Radon-type results are proposed for graphs which are 
analogous to their original versions in discrete convex geometry. For further insight, one can look into the papers \cite{Jamison1984, Duchet1988, Duchet1987, Chvatal1975}. However, this standard set-theoretic definition restricts many further explorations. Motivated by discrete function theory, we impose weights on the vertices and propose several new convex terminologies for graphs. These newly introduced definitions help us to obtain Krein-Milman-type theorems (see \cite{Krein}), Ulam-type stability results, allow us to demonstrate the underlying relationship between sequential convexity and trees of a specific cardinality, and also enhance our understanding of weight allocations to vertices.\\ 

Mitrinovi\'c first introduced the terminology of convex sequence or sequential convexity
in his book \cite{Mitrinovicc}.
A sequence $\big(u_i\big)_{i=0}^{\infty}$ is called \textit{convex} if it satisfies the following discrete functional inequality
\Eq{905}{
2u_i\leq u_{i-1}+u_{i+1}\qquad \mbox{for all}\qquad i\in\N.
}
Arithmetic, geometric, Fibonacci, partition, factorial, and many other well-known classes of sequences satisfy the above inequality. Analysis of generalised, higher-order, and approximate versions of sequential convexity has been conducted over the last few decades. Some fundamental studies related to sequential convexity can be found in the articles \cite{Essen, pecaric, GauSte, Debnath, Latreuch} and their references. Motivated by this, we formulate the following versions of convexity that can be implemented on graphs.\\

In the graph $G$, the sets $V$ and $E$ denote the collections of vertices and edges, respectively. We assume $w:V\to\R_+$ is a non-negative function. A weighted graph $G(V,E,w)$ is said to be \textit{$n$-convex} if for every vertex $v$ with at least $n$ number of adjacent vertices $v_1,\cdots,v_n$, satisfy the following discrete functional inequality
\Eq{222}{
n\cdot w(v)\leq w(v_1)+\cdots+w(v_n),
}
where $w(v), w(v_1),\cdots,w(v_n)$ represent the respective vertex weights. If the inequality holds with the $<$ sign, then we call it a strictly \textit{$n$-convex graph}. \\

We go through several structural characteristics of graphs that exhibit $n$-convexity. We demonstrate that an $n$-regular or a complete graph $K_{n+1}$ can not be strictly $n$-convex. We prove that if a weighted graph is strictly $n$-convex, then the maximum value is always attained at a vertex with degree at most $n-1$. Additionally, we demonstrate that if a graph is strictly $ 2$-convex, then it must be a tree. \\

In our research, one interesting finding is the implementation of Ulam-type stability results in a graph setting. In function theory, Ulam-type stability describes the phenomenon that approximate solutions of a functional equation (or inequality) remain close to exact solutions, demonstrating robustness of the equation(inequality) under small perturbations. Some of the classical work of Hyers and Ulam can be found in the papers \cite{Ulam,Hyers}. For the past 75 years, mathematicians have been proposing new versions of stability results by analysing functional, ordinary, and partial differential equations and inequalities for various classes of real functions. Motivated by these, we investigate Ulam-type stability for graphs. Under some minimal assumptions, we show that  for a fixed $\varepsilon>0$, if every vertices $v$ of the graph $G(V,E,w)$ with degree at least $n$ satisfy the following discrete functional inequality
\Eq{*}{
 w(v)\leq \dfrac{w(v_1)+\cdots+w(v_n)}{n}+\varepsilon \qquad \qquad (v_{1},\cdots,v_n \mbox{ are arbitrary adjacent vertices of $v$}),
}
then there exists an $n$-convex graph $G(V,E,\widetilde{w})$ such that $\|w-\widetilde{w}\|_{\infty}\leq\varepsilon/2$ holds. For readability purposes, we break the proof into several small parts and propose some minorant and sandwich-type results.
\\

We explicitly investigate the case $n=2$, i.e., $2$-convex trees. A weighted tree $T(V,E,w)$ is labelled as \textit{type-$1$ convex} if for any internal vertex $v\in T(V)$ with any two adjacent vertices $v_1$ and $v_2$ satisfies the following discrete inequality 
\Eq{77}{
2w(v)\leq w(v_1)+w(v_2).
}
If the above inequality holds with `$<$' sign, we call $T$ possesses \textit{type-$1$ strict convexity}.\\

We show that the maximum weight of such a tree is always carried by one of its leaf vertices. We prove that, under minimal assumptions, the root vertex of a rooted tree attains the minimum weight. Additional investigations have also been performed to explore other structural properties.\\

On the other hand, a rooted tree $T(V,E,w)$ is called \textit{type-$2$ convex} if every internal vertex $v$ that possesses at least two child nodes $v_1$ and $v_2$ satisfies the inequality \eq{77}. Although the definition applies to arbitrary rooted trees, it is particularly natural in the case of binary trees, where each vertex has at most two child nodes. This setting is relevant to processes involving successive binary splitting or decomposition. Examples include the division of a quantity into two parts, valuation problems arising from corporate de-mergers, and certain branching processes that occur in biological models. Such situations motivate the study of type-$2$ convexity on binary trees. Also, this newly introduced convexity on trees has a deep connection with sequential convexity.\\

Let $T(V,E)$ be a non-weighted rooted tree with $k$ levels (excluding the root) and each level has $n_{i}$ number of vertices $(i\in\{1,\cdots,k\})$ and $\big(u_n\big)_{n\in\N}$ be a monotone (increasing) sequence with $n_1+\cdots+n_k+1$ elements. Then we can assign weights to the $V(T)$ as follows\\

\begin{enumerate}
\item Assign the first element of $\big(u_{n}\big)_{n\in\N}$ as weight to the root vertex of $T$\\

\item Assign the next $n_1$ elements of $\big(u_{n}\big)_{n\in\N}$ as weights to the vertices of the first level from left to right in order.\\

\item Continue this assignment level-wise until the allocation is complete.\\
\end{enumerate}

We observe that for any internal vertex $v$ at level $i$, its child nodes $v_1$ and $v_2$ lie at level $i+1$. Hence, from the monotonicity of the sequence, we have
\Eq{*}{
w(v)\leq w(v_1)\qquad \mbox{and} \qquad w(v)\leq w(v_2).
}

This ensures that \eq{77} holds, thereby turning $T(V,E)$ into a type-$2$ weighted convex tree. In fact, for any given set of numbers, we can first sort them in ascending order and then apply the above procedure to obtain the desired convexity-preserving allocation. Thus, this methodology can also be treated as an exact algorithm. However, for very large numbers, this procedure is computationally infeasible.\\

We demonstrate that if the sequence $\big(u_i\big)_{i=1}^{2^n-1}$ is convex, then it can be embedded into an unweighted perfect binary tree $T(V,E)$ with $n-1$ levels in such a way that it can be turned into a type-$2$ convex tree. To achieve this, we define a bijection $w:V(T)\to U$ without making any alteration to the sequence. This can also be seen as a combinatorial optimisation problem in which the goal is to assign resources to achieve an optimal outcome while minimising allocation costs. Such assignment problems are broadly studied in transportation analysis.\\

Towards the end, we show that for a fixed $\varepsilon>0$, if an arbitrarily chosen internal vertex $v\in V(T)$ together with its any arbitrary descendant vertex $v'$ equipped with the functional inequality $w(v)\leq w(v')+\varepsilon$, then there exists a type-$2$ convex tree $T(V,E,\widehat{w})$ such that $\|w-\widehat{w}\|_{\infty}\leq\varepsilon/2$ holds.\\

We start our investigations with $n$-convexity on graphs.

\section{On $n$-convexity and type-$1$ convexity}
Throughout this section, we assume that all graphs are simple and connected, although this assumption may be relaxed in some cases. The symbol $\Delta(G)$ is used to denote the maximum degree of the graph $G(V,E,w)$. 
We begin with the following inclusion property.
\Prp{889}{
If $\Delta(G)\geq n+1$ and possesses $n$-convexity, then $G(V,E,w)$ is $(n+1)$-convex.
}
\begin{proof}
To prove the result, we assume $v\in V$ is an arbitrary vertex with at least $n+1$ adjacent vertices. Let $v_1,\cdots,v_n,v_{n+1}$ are some of those vertices. Since $G(V,E,w)$ is $n$-convex, it will satisfy the inequality \eq{222}. Without loss of generality, we can assume that
$w(v_1)\leq \cdots\leq w(v_n) \leq w(v_{n+1}).$ Then the $n$-convexity of $G$ implies $w(v)\leq w(v_n).$ If not, we have the following system of inequalities
\Eq{*}{
w(v_1)<w(v),\cdots\cdots\mbox{and}\,\,, w(v_n)<w(v)\qquad \mbox{which implies}\qquad \sum_{i=1}^{n}w(v_i)<n\cdot w(v).
}
This results in a contradiction. Hence, $w(v)\leq w(v_n)\leq w(v_{n+1})$ is obvious.
This together with \eq{222}, we can compute the following inequality
\Eq{*}{
(n+1)\cdot.w(v)=n\cdot w(v)+w(v)\leq w(v_1)+\cdots+w(v_n)+w(v_{n+1}).
}
Since $v\in V$ is arbitrary, we can conclude that $G$ possesses $(n+1)$-convexity. This completes the proof.
\end{proof}
The next result demonstrates that if $H\subset G$ is $k$-regular, then $G(V,E,w)$
cannot be a strict $n$-convex graph for any $n\leq k$.
\Prp{890}{Let $G(V,E,w)$ be a weighted be $k$-regular graph. Then $G$ cannot possess strict $n$-convexity for any $n\leq k$.
}
\begin{proof}
Let the graph $G(V,E,w)$ has $m$ vertices. If possible, we assume that there exists an $n\leq k$ such that $G$ possesses $n$-convexity. Then by the above establishment (\prp{889}), $G(V,E,w)$ is $k$-convex as well. Let 
$v_i \in V$ be an arbitrary node and $v_{_{i}}^{^{1}},\cdots,v_{_{i}}^{^{k}}$ are $k$ adjacent vertices of it. Then, due to strict $n$-convexity, it satisfies the following inequality
\Eq{*}{
k\cdot w(v_i)<w(v_{_{i}}^{^{1}})+\cdots+w(v_{{_i}}^{^{k}}).
}
Considering all such $v_i^{'}s$, we obtain a system of $m$ inequalities. Since each $v_i\in V$ is adjacent to exactly $k$ vertices, summing up all such inequalities side by side, we arrive at
\Eq{*}{
k\sum_{i=1}^{m}w(v_i)<k \sum_{i=1}^{m}w(v_i).
}
The inequality above is contradictory; hence, our assumption is wrong. This establishes the result.
\end{proof}
The next theorem can also be treated as a Krein-Milman-type result for graph theory.
\Prp{0}{If $G(V,E,w)$ is strictly $n$-convex, then the maximum value is attained by a vertex which has a degree less than $n$.}
\begin{proof}
To establish the proposition, we assume $v_{0}\in V$ be a vertex such that 
$\underset{v\in V}{\max}\Big({w(v)}\Big)=w(v_0)$. If possible, let degree of $v_0$ be atleast $n$. In other words, there exists at least $n$ adjacent vertices $v_1,\cdots,v_n$ of $v_0$ that satisfy the following inequality
\Eq{*}{
n\cdot w(v_0)< w(v_1)+\cdots +w(v_n).
} 
The non-negativity of weights and strict $n$-convexity ensure that there exists a vertex $v_i$ such that $w(v_0)<w(v_i)$ holds. This yields a contradiction and proves our assertion.
\end{proof} 
However, the minimum weight can also lie in a vertex of degree $n$ or more. For illustration, we have the following graph.
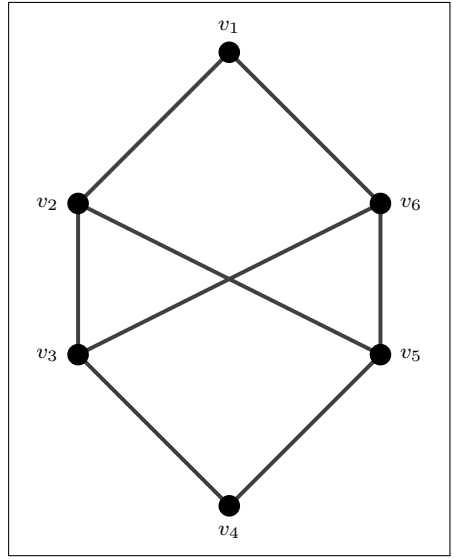
\begin{center}
\fbox{
\begin{tikzpicture}
\Vertex[size=0.25, color=black, label=$v_1$, position=above]{A}
\Vertex[x=-2,y=-2, size=0.25, color=black, label=$v_2$,  position=left]{B}
\Vertex[x=-2, y=-4, size=0.25, color=black, label=$v_3$,  position=left]{C}
\Vertex[x=0, y=-6,  size=0.25, color=black, label=$v_4$,  position=below]{D}
\Vertex[x=2, y=-4, size=0.25, color=black, label=$v_5$,  position=right]{E}
\Vertex[x=2, y=-2, size=0.25, color=black, label=$v_6$,  position=right]{F}
\Edge(A)(B)
\Edge(B)(C)
\Edge(C)(D)
\Edge(D)(E)
\Edge(E)(F)
\Edge(F)(A)
\Edge(C)(F)
\Edge(E)(B)
\end{tikzpicture}
}
\end{center}
\begin{center}
$w(v_1)=10$, $w(v_2)=1$, $w(v_3)=2$, $w(v_4)=20$, $w(v_5)=5$, $w(v_6)=3$\\
\end{center}
\hspace{1cm}

With the assigned weights, the graph above is $3$-convex. The maximum weight is carried by $v_4$, a vertex of degree $2$ as established in \prp{0}. In contrast, the minimum weight is borne by $v_2$, a vertex with degree $3$. This validates our statement.\\

The following result follows directly from \prp{890} and \prp{0}.
\Cor{0}{The complete graph $K_{n+1}$ does not possess strict $n$-convexity.}

Before stating the next results, we first need to introduce several notations. For the graph $G(V,E,w)$, the symbol $V_{{n}}$ will be used to denote the set of vertices having degree at least $n$. The symbol $\A(v)$ denotes the set of all vertices adjacent to the vertex $v$. We now present the Ulam-type stability result. For clarity, we divide the proof into several steps. First, we establish a result concerning convex minorants. Next, we prove a sandwich-type result. Finally, combining these two results, we derive the Ulam-type stability theorem as a consequence.\\

The graph $G(V,E,w_{0})$ is called a \textit{minorant} of $G(V,E,w)$ if the inequality
$w_0(v)\leq w(v)$ holds for all $v\in V$. The study of minorants is crucial in approximation theory. In particular, non-trivial (non-constant) minorants that satisfy properties such as monotonicity, convexity, or subadditivity, etc. are especially interesting.
\Prp{624}{
Let $G(V,E,w)$ be an arbitrary graph with $\Delta(G)\geq n$. Then there exists an $n$-convex minorant $G(V,E,\widetilde{w})$ (possibly non-trivial) of $G(V,E,w)$.}
\begin{proof}
From our assumptions on $G$, it is clear that $V_n$ is non-empty.
To prove the assertion, first we define the iteration $w_k$ $(k\in\N)$ on the weight function $w$ as follows
\begin{small}
\Eq{329}{
w_{_0}(v):=w(v) \quad \mbox{and}\quad 
w_{_{k}}(v):=\begin{cases}
\qquad \qquad \quad \quad \quad \quad w(v)\,\,\, &\mbox{if}\,\,\, v\in V\setminus{V_{n}}\\
\min\bigg\{w_{_{k-1}}(v)\,\,,\,\,\underset{u_i\in\A(v)}{\inf}\dfrac{w_{_{k-1}}(u_1)+\cdots +w_{_{k-1}}(u_n)}{n}\bigg\}\,\,\, &\mbox{if}\,\,\, v\in {V_{n}}
\end{cases}
}
\end{small}
From the construction, it is clear that for any $v\in V$, the sequence 
$\big(w_k(v)\big)_{k\in\N}$ is decreasing and bounded below by 0. This ensures the convergence of the sequence. Using this, we define the function $\widetilde{w}:V\to\R_+$ as 
$\widetilde{w}(v)=\underset{{k\to\infty}}{\lim}w_{k}(v)$
and claim that $\widetilde{w}$ is a convex minorant of $w$.\\

By definition, for any $v\in V$ the inequality $\widetilde{w}(v)\leq w(v)$ is obvious. To prove $n$-convexity, we consider two cases. For $v\in V\setminus V_{n}$, there is nothing to show. On the other hand, if $v\in V_{n}$, we assume $n$ vertices $v_1,\cdots,v_n\in \A(v)$ arbitrarily. Then we can compute the following inequalities
\Eq{*}{
\widetilde{w}(v)=\lim_{k\to\infty}w_k(v)&=\lim_{k\to\infty}\Bigg(\min\bigg\{w_{_{k-1}}(v)\,\,,\,\,\underset{u_i\in\A(v)}{\inf}\dfrac{w_{_{k-1}}(u_1)+\cdots +w_{_{k-1}}(u_n)}{n}\bigg\}\Bigg)\\
&\leq \lim_{k\to\infty}\Bigg(\min\bigg\{w_{_{k-1}}(v)\,\,,\,\,\dfrac{w_{_{k-1}}(v_1)+\cdots +w_{_{k-1}}(v_n)}{n}\bigg\}\Bigg)\\
&\leq \lim_{k\to\infty}\,\,\,\dfrac{w_{_{k-1}}(v_1)+\cdots +w_{_{k-1}}(v_n)}{n}\\
&=\dfrac{\widetilde{w}(v_1)+\cdots+\widetilde{w}(v_n)}{n}
}
This validates the convexity of $\widetilde{w}:V\to\R_+$ and establishes that $G(V,E,\widetilde{w})$ is an $n$ convex minorant of $G(V,E,{w})$.
\end{proof}
Using the above proposition, we can obtain the following sandwich-type result for graphs.
\Cor{6231}{Let $G(V,E,w')$ and $G(V,E,w)$ be two weighted graphs  with $\Delta(G)\geq n$ such that $w'(v)\leq w(v)$ holds for all $v\in V$. Additionally, the for all $v\in V_n$, the following condition applies
\Eq{330}{
 w'(v)\leq \dfrac{w_k(u_1)+\cdots w_k(u_n)}{n} \quad \mbox{for all} \quad k\in\N\cup\{0\}\quad \Big(u_1,\cdots,u_n\in \A(v) \Big),
}
where $w_k$ is the iteration defined in \eq{329}. Then there exists an $n$-convex graph $G(V,E,\widetilde{w})$ that satisfies the inequality $w'(v)\leq \widetilde{w}(v)\leq w(v)$ for all $v\in V$.}
\begin{proof}
To prove the corollary, we choose any $v\in V$ and then consider the decreasing sequence $\Big(w_k(v)\Big)_{k\in\N}$ as in \eq{329}, which is bounded below by $w'(v)$. Hence, the two conditions of \eq{330} ensures that the weight function $\widetilde{w}: V\to\R_+$ defined as $\widetilde{w}(v):=\underset{{k\to\infty}}{\lim}w_{k}(v)$ leads to the $n$-convex graph $G(V,E,\widetilde{w})$ satisfying $w'(v)\leq \widetilde{w}(v)\leq w(v)$ for all $v\in V$. This completes the proof.
\end{proof}
Let $\varepsilon>0$ be fixed. A weighted graph $G(V,E,w)$ with $\Delta(G)\geq n$ is said to be \textit{$(n,\varepsilon)$-convex or approximately convex} if for each vertex $v\in V_n$ with any of its $n$ adjacent vertices $v_1,\cdots,v_n$ satisfy the following discrete functional inequality
\Eq{*}{
w(v)\leq\dfrac{w(v_1)+\cdots+w(v_n)}{n}+\varepsilon.
} 
Finally, we are now able to propose a Ulam-type stability result for $n$-convex graphs. 
\Thm{959}{Let $\varepsilon>0$ be fixed. If $G(V,E,{w})$ is a $(n,\varepsilon)$-convex graph such that it satisfies the following discrete functional inequality
\Eq{8703}{
 w(v)\leq \dfrac{w_k(u_1)+\cdots w_k(u_n)}{n}+\varepsilon \quad \mbox{for all} \quad k\in\N\cup\{0\}\quad \Big(u_1,\cdots,u_n\in \A(v) \Big),
}where $w_k$ is the iteration defined in \eq{329},
then there exists an $n$-convex graph $G(V,E,\widetilde{w})$ such that the following norm inequality
$\|w-\widetilde{w}\|_{\infty}\leq\dfrac{\varepsilon}{2}$ holds. Conversely, if $G(V,E,\widetilde{w})$ is a convex graph that satisfies the inequality $\|w-\widetilde{w}\|_{\infty}\leq\dfrac{\varepsilon}{2}$, then $G(V,E,w)$ is a $(n,\varepsilon)$-convex graph.
}
\begin{proof}
First, we observe that the inequality \eq{8703} can also be represented as follows
\begin{small}
\Eq{8704}{\max\big\{{w}(v)-{\varepsilon}/{2}\,,\,0\,\}\leq \dfrac{\Big({w_k}(u_1)+\varepsilon/2\Big)+\cdots+\Big({w_k}(u_n)+\varepsilon/2\Big)}{n}\quad \mbox{for all} \quad k\in\N\cup\{0\}\quad \bigg( u_i\in\A(v)\bigg)
}
\end{small}
We now introduce the weight function $w':=\max\{{w}-\varepsilon/2\,,\,0\}$ and replace  
$ w_k:={w_k}+\varepsilon/2$.
This together with \eq{8704} yields \eq{330}. Hence, utilising \prp{624}, and \cor{6231}, we conclude the existence of an $n$-convex graph $G(V,E,\widetilde{w})$ such that for all
$v\in V$, the following inequality holds
\Eq{*}{
w'(v)=\max\Big\{{w}(v)-\varepsilon/2\,,\,0\Big\}\leq\widetilde{w}(v)\leq{w}(v)+\varepsilon/2.\
}
From the above inequality, we conclude the following 
\Eq{2026}{
{w}(v)-\varepsilon/2\leq\widetilde{w}(v)\leq{w}(v)+\varepsilon/2 \qquad \mbox{for all}\quad v\in V.
}
This validates the norm inequality 
$\|w-\widetilde{w}\|_{\infty}\leq\dfrac{\varepsilon}{2}$ and  proves the first-part of the assertion.\\

To establish the second part of the theorem, let $G(V,E,\widetilde{w})$ be an $n$-convex graph that satisfies the inequality \eq{2026}. Let $v\in V_n$ be arbitrary and $v_1,\cdots,v_n\in\A(v)$. Now, first using the left-most inequality of \eq{2026} and then the right-most inequality of it, we can compute the following
\Eq{*}{
w(v)\leq \widetilde{w}(v)+\dfrac{\varepsilon}{2}
&\leq \dfrac{\widetilde{w}(v_1)+\cdots+\widetilde{w}(v_n)}{n}+\dfrac{\varepsilon}{2}\\
&\leq \dfrac{\big({w}(v_1)+\varepsilon/2\big)+\cdots+\big({w}(v_n)+\varepsilon/2\big)}{n}+\dfrac{\varepsilon}{2}\\
&=\dfrac{w(v_1)+\cdots+w(v_n)}{n}+\varepsilon.
}
This shows that $G(V,E,w)$ is $(n,\varepsilon,)$-convex and completes the proof.
\end{proof}
From now on, we will only focus on studying weighted convexities on trees.
\Prp{1}{A strictly $2$-convex graph is a tree.}
\begin{proof}
Let $G(V,E,w)$ be a strictly $2$-convex graph. Suppose, if possible, $G$ is not a tree. Without loss of generality, we assume that there exists a cycle, $C_n=(v_1,v_2\cdots, v_n,v_1)$ of $n$ vertices in $ G$. Now, the strict convexity of $G$ yields,
\Eq{*}{
2w(v_2)<w(v_1)+w(v_3).
}
We claim that, $w(v_2)<\max\{w(v_1),w(v_3)\}. $ If not, then we have
\Eq{*}{
w(v_1)+w(v_3)\leq 2 \max\{w(v_1),w(v_3)\}\leq 2w(v_2).
}
This contradicts the strict $2$-convexity of $G$. Therefore, without loss of generality, we can assume that
\Eq{39}{
w(v_1)\leq w(v_2)<w(v_3).
} 
Again, due to strict $2$-convexity, the vertices $v_2, v_3$ and $v_4$ satisfy the following inequality
\Eq{*}{
2w(v_3)<w(v_2)+w(v_4).
}
And analogously together with the inequality \eq{39}, we can also conclude
\Eq{*}{
w(v_2)<w(v_3)<w(v_4).
}
Proceeding this way, we will eventually obtain
\Eq{*}{
w(v_1)\leq w(v_2)<w(v_3)<w(v_4)<\cdots<w(v_n)<w(v_1).
}
This results in a contradiction. Instead of the initial assumption of \eq{39}, by starting with other possibilities, we can present the same irregularity. Hence, $G(V,E)$ is a simple, connected graph with no cycles. In other words, $G(V,E,w)$ is a tree.
The claim follows.
\end{proof}
The next proposition demonstrates that the vertices attaining the maximum and minimum values of a type-$1$ convex tree are predictable. Hence, the finding can be treated as a Krein-Milman-type theorem for trees.
\Prp{2}{ In a type-$1$ strictly convex tree, the maximum weight is carried by at least one of the leaf vertices. Additionally, besides type-$1$ strict convexity of the tree, if the tree is rooted and the weight of the root vertex is not greater than the weights of any of its adjacent vertices, then the root vertex attains the minimum weight. 
}
\begin{proof}
Let $T(V,E,w)$ be a type-$1$ strictly convex tree. Since a type-$1$ strictly convex tree refers to a $2$-convex graph, by \prp{0}, the maximum weight is carried by a vertex of degree $1$. That is, a leaf vertex bears the maximum weight. This establishes the first assertion.\\

To prove the second assertion, let $V_{\ell}\subset{V}$ be the set of leaf vertices. Let $v_0\in V_{\ell}$ be arbitrary and $\tilde{v}\in V$ be the root of the tree. If all the adjacent vertices of $\tilde{v}$ are leaves, we have nothing to prove. Otherwise, we consider an arbitrary path $P$ from $v_0$ to $\tilde{v}$  as follows
\Eq{*}{P:=(v_0,v_1,,\cdots,v_{n-1},v_n,v_{n+1}(=\tilde{v})) \quad \mbox{where}\quad v_1,\cdots,v_n \quad \mbox{are the internal nodes of $P$}.}

Then, using the assumptions on the root node $\tilde{v}$, we have the following two possibilities
\Eq{413}{
 w(\tilde{v})=w(v_{n+1})\leq w(v_{n-1})\leq w(v_{n})\quad \mbox{and}\quad 
 w(\tilde{v})=w(v_{n+1})\leq w(v_n)\leq w(v_{n-1}).
}
From the first inequality of \eq{413}, we have
$
w(v_{n-1})+w(v_{n+1})\leq 2w(v_n),
$
which is a contradiction to the type-$1$ strict convexity of $T$. Hence, we need to analyse only the second inequality of \eq{413}.\\

Using the convexity of $T(V,E,w)$, we can obtain the following discrete inequality 
\Eq{*}{
w(\tilde{v})= w(v_{n+1})\leq w(v_{n})<\cdots <w(v_0).
}
Since $v_0$ and $P$ are arbitrarily  chosen, the weighted monotonicity of the path $P$ establishes the following inequalities
\Eq{*}{
\max_{v\in V}{w(v)}=\max_{v\in V_\ell}{w(v)}\qquad \mbox{and} \qquad
\min_{v\in V}{w(v)}={w(\tilde{v})}.
}
This validates the second assertion and completes the proof.
\end{proof}
In the next section, we will establish a relationship between sequences and graphs, the two primary branches of discrete mathematics.
\section{On type-$2$ convexity}
Before proceeding, we recall the definitions of proper and perfect binary trees. 
A binary tree is said to be \textit{proper} if each node has either zero or two children. On the other hand, a \textit{perfect binary tree} is a special type of proper binary tree that satisfies both of the following conditions
\begin{itemize}
\item All internal nodes have exactly two children.
\item All leaf nodes are at the same level. In other words, every root to leaf path has the same length.
\end{itemize}

\Thm{09}{
Let $\big(u_i\big)_{i=1}^{2^n-1}$ $(n\in\N)$ be a convex sequence. Then it can be embedded into an unweighted perfect binary tree of $n-1$ levels to make it type-$2$ convex.
}
\begin{proof}
The proof of the theorem contains several steps. First, we need to characterise convex sequences. The characterisation can also be found in \cite{Gosw}. But for clarity, we decide to elaborate on it. Then we represent a unique way to express any number $k\in[1,2^{n}-1]$. Finally, we define a mapping $w:V\to\big(u_i \big)_{i=1}^{2^n-1}$ such that the $T(V,E)$ turns into a type-$2$ convex tree.\\

\textbf{The First Step: We validate the following statement}\\

``A sequence $\big(u_i\big)_{i=1}^{\infty}$ is convex if and only if for any $p,q,r,s\in\N$ with $p<q\leq r<s$, it satisfies the following discrete functional inequality
\Eq{1111}{
u_q+u_r\leq u_p+u_s \qquad \mbox{provided}\qquad q+r=p+s\,\,.\,^{^"}
}
At first, we assume that the sequence $\big(u_i\big)_{i=1}^{\infty}$ is convex and $p,s\in\N$ with $p<s$ are fixed. This leads us to the following system of inequalities
\begin{equation}
\tag{P}
{ u_{p+1}-u_p\leq u_{p+2}-u_{p+1}}
\end{equation}
\begin{equation}
\tag{P+1}
{ u_{p+2}-u_{p+1}\leq u_{p+3}-u_{p+2}}
\end{equation}
$$\vdots \qquad \qquad \vdots$$
\begin{equation}
\tag{S-3}
{ u_{s-2}-u_{s-3}\leq u_{s-1}-u_{s-2}}
\end{equation}
\begin{equation}
\tag{S-2}
{ u_{s-1}-u_{s-2}\leq u_{s}-u_{s-1}}.
\end{equation}
Adding up all the inequalities, side by side, we obtain 
\Eq{7788}{
u_{p+1}+u_{s-1}\leq u_p+u_s.
} 
Now, excluding the inequalities $P$ and $S-2$ from the above system and then adding up all inequalities side by side, we arrive at
\Eq{7789}{
u_{p+2}+u_{s-2}\leq u_{p+1}+u_{s-1}.
}
Similarly, excluding the inequalities $P$, $P+1$, $S-3$, and $S-2$, we sum up all the remaining inequalities and obtain
\Eq{7790}{
u_{p+3}+u_{s-3}\leq u_{p+2}+u_{s-2}.
}
We keep continuing the process, and each time we find new inequalities similar to \eq{7788}, \eq{7789}, and \eq{7790}. Depending upon the values of $p$ and $q$, we will end up in a system of inequalities that can be summarised as follows
\Eq{*}{
2u_{\frac{p+s}{2}}\leq u_{\frac{p+s}{2}-1}+ u_{\frac{p+s}{2}+1}\leq &\cdots \leq u_{p+1}+u_{s-1}\leq u_{p}+u_{s} \quad\mbox{if $p+s$ is even}\\
&\mbox{or}\\
u_{\frac{p+s-1}{2}}+ u_{\frac{p+s+1}{2}}\leq u_{\frac{p+s-1}{2}-1}+ u_{\frac{p+s-1}{2}+1}\leq&\cdots \leq u_{p+1}+u_{s-1}\leq u_{p}+u_{s} \quad\mbox{if $p+s$ is odd.}
}
The above inequalities show all the 
possible combinations of $q,r\in\N$ such that if $p<q\leq r<s$ holds, then \eq{1111} is satisfied.\\

Conversely if \eq{1111} holds for all $p,q,r,s\in\N$, then by choosing $p=i-1$, $q=r=i$ and $s=i+1$, ($i\in \N\cap[2,\infty[\,$); we get the inequality \eq{905}.
This yields convexity of the sequence $\big(u_i\big)_{i=1}^{\infty}$ and establishes the characterisation. \\

\textbf{The second Step: We prove the following unique representation}\\

"For any $k\in\N\cap[1, 2^n-1]$, there is a unique representation of $k$, given by
\Eq{71}{
k=m\cdot 2^{n-i-1}\quad\mbox{where}\quad i\in[0,n-1]\cap(\N\cup\{0\}) \quad\mbox{and} \quad m\in \N\cap[1, 2^i-1]\quad\mbox{is odd}\,\,.
^{{"}}
}
For any $k\in\N\cap[1, 2^n-1]$, primarily there are two possibilities. If $k$ is an odd number, then we have the following representation
\Eq{1112}{
k=k\cdot 2^{0}=k\cdot 2^{(n-1)-(n-1)}
}
which align with the claimed representation of \eq{71} with $m=k$ and $i=n-1$.\\

On the other hand, any even $k\in\N$ can be expressed in the following form
\Eq{72}{
k=m\cdot 2^{j}\quad \mbox{with} \quad m, j\in \N \quad (m\,\,\,\, \mbox{is odd}).
}
If $k\in [1,2^n-1]$; then $j\geq n$ implies the following contradiction.
$$k=m\cdot 2^{j}\geq m\cdot 2^{n}\geq 2^{n};$$
Hence, using this in \eq{72}, we can improve the representation of any even $k\in \N\cap[1,2^n-1]$ as follows
\Eq{*}{
k=m\cdot 2^{n-i-1}\quad \mbox{where}\quad i\in[0,n-2]\cap(\N\cup\{0\})\quad \mbox{and} \quad m\in \N \quad \mbox{is odd}.
}

Again, under the above mentioned restrictions on $k$, if $m> 2^{i+1}$, then $k=m\cdot 2^{n-i-1}> 2^{i+1}\cdot 2^{n-i-1}=2^n$, This also results in a contradiction. Hence, for any even $k\in \N\cap[1,2^n-1]$, we have the following refined expression

\Eq{1113}{
k=m\cdot 2^{n-i-1}\quad \mbox{where} \quad i\in[0,n-2]\cap(\N\cup\{0\})\quad\mbox{and}\quad m\in\N\cap[1,2^{i+1}-1] \quad \mbox{is odd}.
}
The combining case-wise representations in \eq{1112} and \eq{1113} yields the validity of \eq{71}.\\

Now, we show the uniqueness of this representation. If possible, let for any $k\in\N\cap[1, 2^n-1]$, there are two distinct representations as follows
$$k=m_12^{n-i-1}\qquad \mbox{and}\qquad k=m_22^{n-j-1}$$
For $i=j$ or $m_1=m_2$, the contradiction is evident. Hence, there is nothing to show. For the case, $ m_1\neq m_2$ and
$i\neq j$, we can compute the following equality
$$m_12^{n-i-1}=m_22^{n-1-j} \quad \mbox{or}\quad \dfrac{m_1}{m_2}=2^{i-j}.$$
This yields that at least one of $m_1$ and $m_2$ is even, and this results in a contradiction. Hence, the representation in \eq{71} is unique.\\

For simplicity, one can also look into representation of $k\in\N\cap[1,2^n-1]$ of \eq{71} in the following format as-well
\Eq{7}{
k=(2l-1)\cdot 2^{n-i-1}\quad\mbox{where}\quad i\in[0,n-1]\cap(\N\cup\{0\}) \quad \mbox{and} \quad l\in \N\cap[1, 2^{i-1}].
}\\
\vspace{5mm}
\textbf{The Third Step: We define the mapping between the two spaces}

Finally, we assign weights to each vertex of the unweighted perfect binary tree $T(V,E)$. Let $T$ has $n-1$ levels (excluding the root). Since $T$ is perfect, each level has exactly $2^{i}$ $(i\in\{0,1,\cdots,n-1\})$ nodes totalling $2^{n}-1$ elements. The symbol $v_l^i$ is used to denote the $l^{th}$ vertex (from left) in the $i^{th}$ level. For any internal node $v_l^i$ of $T(V,E)$, the child nodes will be $v_{_{2l-1}}^{^{i+1}}$ and  $v_{_{2l}}^{^{i+1}}$ lying in the $(i+1)^{th}$ level.\\

From \eq{7}, the sequence $\big(u_i\big)_{i=1}^{2^n-1}$  is also representable as follows
$$\big(u_{_{1.{2^0}}}, u_{_{{1.{2^1}}}}, u_{_{{3.2^0}}}, u_{_{1.{2^2}}},\cdots, u_{_{{2^n-1}.2^0}}\big).$$

Now we define a mapping $w:V(T)\to \big(u_i\big)_{_{i=1}}^{^{2^n-1}}$ as follows
\Eq{47}{
w\Big(v_{_l}^{^i}\Big)=u_{_{(2l-1)\cdot2^{n-i-1}}}\,\, \Big(\mbox{or}\,\, u_{_{m\cdot2^{n-i-1}}}\Big)\quad \mbox{where}\quad l\in \N\cap[1,2^{i-1}]\quad \Big(m\in\N\cap[1,2^i-1]\quad\mbox{is odd}\Big).
}

Hence, for any internal node $v_{_l}^{^i}$, its two child nodes $v_{_{2l-1}}^{^{i+1}}$ and  $v_{_{2l}}^{^{i+1}}$ satisfy the following equalities 

\Eq{*}{
w\Big(v_{_{2l-1}}^{^{i+1}}\Big)=u_{_{(2m-1)\cdot2^{n-i-2}}}\qquad \mbox{and}
\qquad w\Big(v_{_{2l}}^{^{i+1}}\Big)=u_{_{(2m+1)\cdot 2^{n-i-2}}}.
}

Using \eq{1111} together with \eq{47}, we can compute the following inequality
\Eq{*}{
2\cdot w\Big(v_{_l}^{^i}\Big)= 2\cdot u_{_{m\cdot2^{n-i-1}}}
&\leq u_{_{(2m-1)\cdot2^{n-i-2}}}+u_{_{(2m+1)\cdot2^{n-i-2}}}
&=w\Big(v_{_{2l-1}}^{^{i+1}}\Big)+w\Big(v_{_{2l}}^{^{i+1}}\Big).
}
Since $v_{_l}^{^i}$ is an arbitrary internal node, the above inequality confirms type-$2$ convexity in the perfect tree $T(V,E,w)$ and completes the proof.
\end{proof}
The following figure illustrates the above theorem. It shows how the elements of the convex sequence $\big(u_i\big)_{i=1}^{2^n-1}$ are distributed among the nodes of a perfect binary tree $T(V,E)$ with $n-1$ levels such that it becomes a type-$2$ convex tree.
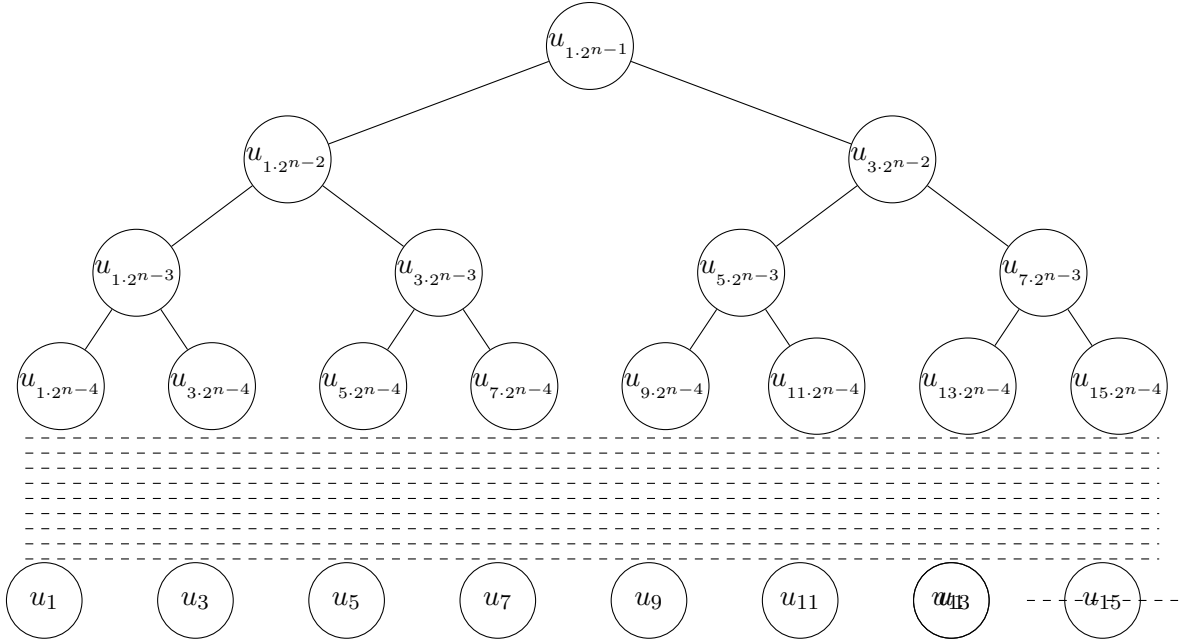
\begin{center}
\begin{tikzpicture}[
  level 1/.style={sibling distance=8cm, level distance=1.5cm},
  level 2/.style={sibling distance=4cm},
  level 3/.style={sibling distance=2cm},
  level 4/.style={sibling distance=1cm},
  every node/.style={circle, draw, minimum size=7mm, inner sep=0pt}
]
\node{$u_{_{1\cdot2^{n-1}}}$}
  child {node{$u_{_{1\cdot2^{n-2}}}$}
    child {node{$u_{_{1\cdot 2^{n-3}}}$}
      child {node{$u_{_{1\cdot 2^{n-4}}}$}}
      child {node{$u_{_{3\cdot 2^{n-4}}}$}}
    }
    child {node{$u_{_{3\cdot 2^{n-3}}}$}
      child {node{$u_{_{5\cdot 2^{n-4}}}$}}
      child {node{$u_{_{7\cdot 2^{n-4}}}$}}
    }
  }
  child {node{$u_{_{3\cdot 2^{n-2}}}$}
    child {node{$u_{_{5\cdot 2^{n-3}}}$}
      child {node{$u_{_{9\cdot 2^{n-4}}}$}}
      child {node{$u_{_{11\cdot 2^{n-4}}}$}}
    }
    child {node{$u_{_{7\cdot 2^{n-3}}}$}
      child {node{$u_{_{13\cdot 2^{n-4}}}$}}
      child {node{$u_{_{15\cdot 2^{n-4}}}$}}
    }
  };
\end{tikzpicture}  
\begin{tikzpicture}
\draw[dashed] (-7.5,0) -- (7.5,0); 
\draw[dashed] (-7.5,0.2) -- (7.5,0.2); 
\draw[dashed] (-7.5,0.4) -- (7.5,0.4); 
\draw[dashed] (-7.5,0.6) -- (7.5,0.6);
\draw[dashed] (-7.5,0.8) -- (7.5,0.8); 
\draw[dashed] (-7.5,1) -- (7.5,1);
\draw[dashed] (-7.5,1.2) -- (7.5,1.2); 
\draw[dashed] (-7.5,1.4) -- (7.5,1.4); 
\draw[dashed] (-7.5,1.6) -- (7.5,1.6);
\end{tikzpicture}
\begin{tikzpicture}[
  every node/.style={circle, draw, minimum size=10mm, inner sep=0pt}
]
  \node (v1) at (0,0) {$u_1$};
  \node (v2) at (2,0) {$u_3$};
  \node (v3) at (4,0) {$u_5$};
  \node (v4) at (6,0) {$u_7$};
  \node (v5) at (8,0) {$u_9$};
  \node (v6) at (10,0) {$u_{11}$};
  \node (v7) at (12,0) {$u_{13}$};
  \node (v8) at (14,0) {$u_{15}$};
  \draw[dashed] (13,0) -- (15,0);

  \node (v9) at (12,0) {$u_1$};
\end{tikzpicture}
\end{center}
The following image illustrates how the convex sequence $\Big(u_i\Big)_{i=1}^{15}$ can be represented as a perfect binary tree.
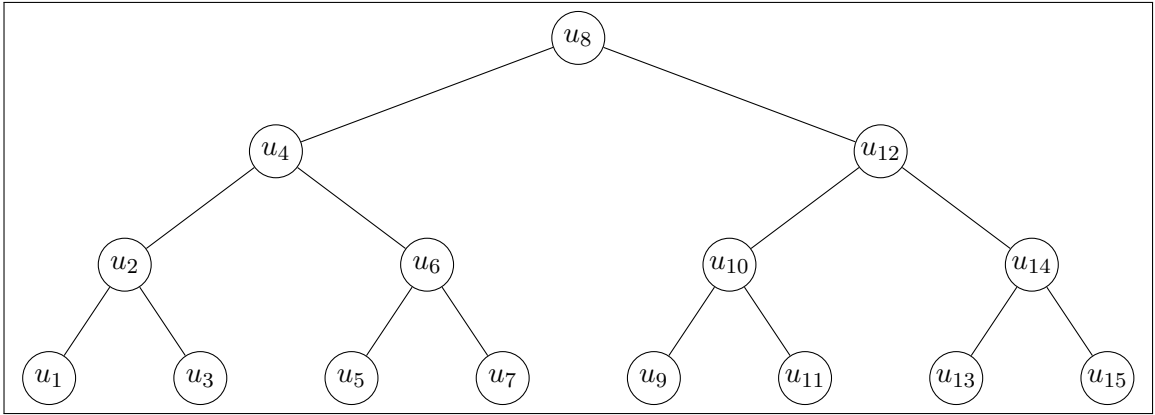
\begin{center}
\fbox{
\begin{tikzpicture}[
  level 1/.style={sibling distance=8cm, level distance=1.5cm},
  level 2/.style={sibling distance=4cm},
  level 3/.style={sibling distance=2cm},
  level 4/.style={sibling distance=1cm},
  every node/.style={circle, draw, minimum size=7mm, inner sep=0pt}
]
\node{$u_8$}
  child {node{$u_4$}
    child {node{$u_2$}
      child {node{$u_1$}}
      child {node{$u_3$}}
    }
    child {node{$u_6$}
      child {node{$u_5$}}
      child {node{$u_7$}}
    }
  }
  child {node{$u_{12}$}
    child {node{$u_{10}$}
      child {node{$u_9$}}
      child {node{$u_{11}$}}
    }
    child {node{$u_{14}$}
      child {node{$u_{13}$}}
      child {node{$u_{15}$}}
    }
  };
\end{tikzpicture}
}
\end{center}
From a type-$2$ convex tree, it is not always possible to derive a convex sequence simply by following the inverse mapping of \eq{47}. Similarly, from a convex sequence, several allocations on an unweighted tree can lead to a type-$2$ convexity. The following examples validate the statement. We consider the sequence $(10,6,3,1,1,4,7)$. Then, as per \thm{09} and the algorithm discussed in the introduction, we obtain the following embeddings.

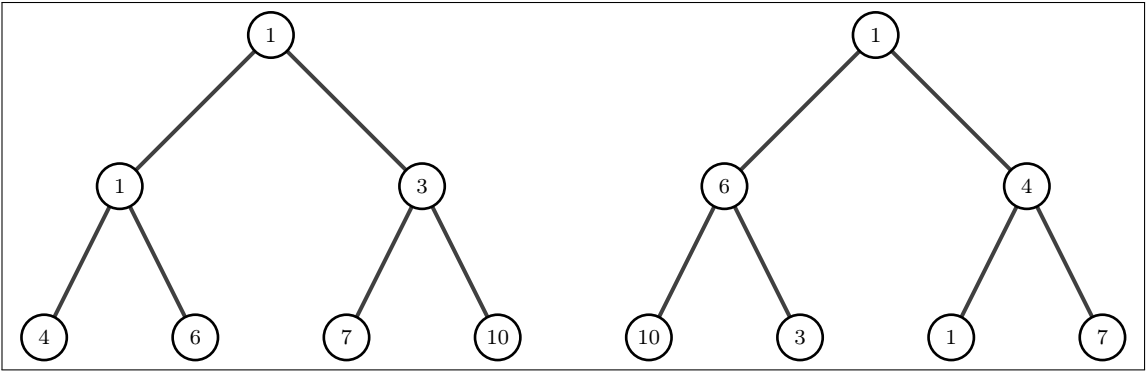
\begin{center}
\fbox{
\begin{tikzpicture}

\Vertex[y=0,x=4, color=white, label=1]{A}
\Vertex[y=-2,x=2, color=white, label=6]{B}
\Vertex[y=-2,x=6, color=white, label=4]{C}
\Vertex[y=-4,x=5, color=white, label=1]{D}
\Vertex[y=-4,x=7, color=white, label=7]{E}
\Vertex[y=-4,x=1, color=white, label=10]{F}
\Vertex[y=-4,x=3, color=white, label=3]{G}
\Edge(A)(B)
\Edge(A)(C)
\Edge(C)(D)
\Edge(C)(E)
\Edge(B)(F)
\Edge(B)(G)

\Vertex[y=0,x=-4, color=white, label=1]{A}
\Vertex[y=-2,x=-6, color=white, label=1]{B}
\Vertex[y=-2,x=-2, color=white, label=3]{C}
\Vertex[y=-4,x=-7, color=white, label=4]{D}
\Vertex[y=-4,x=-5, color=white, label=6]{E}
\Vertex[y=-4,x=-3, color=white, label=7]{F}
\Vertex[y=-4,x=-1, color=white, label=10]{G}
\Edge(A)(B)
\Edge(A)(C)
\Edge(B)(D)
\Edge(B)(E)
\Edge(C)(F)
\Edge(C)(G)

\end{tikzpicture}
}

\end{center}
Besides these two, the following embeddings of the sequence $(10,6,3,1,1,4,7)$ also result in type-$2$ convex graphs.
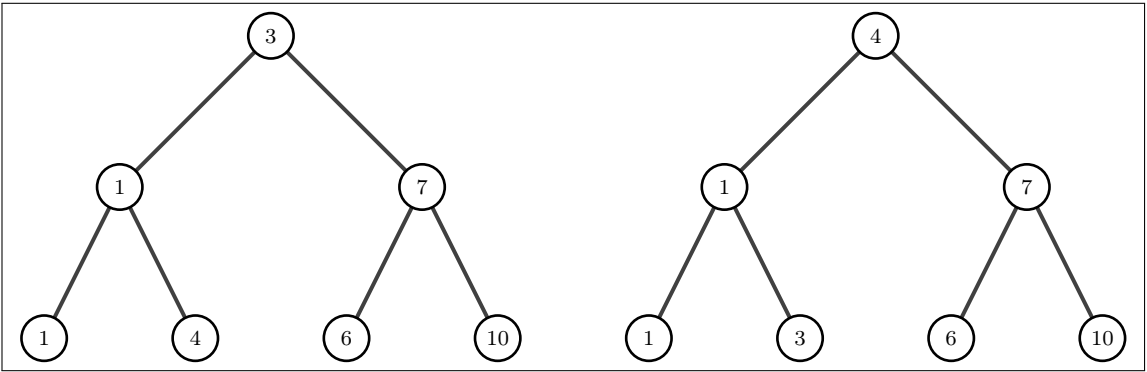
\begin{center}
\fbox{
\begin{tikzpicture}
\Vertex[y=0,x=-4, color=white, label=3]{A}
\Vertex[y=-2,x=-6, color=white, label=1]{B}
\Vertex[y=-2,x=-2, color=white, label=7]{C}
\Vertex[y=-4,x=-7, color=white, label=1]{D}
\Vertex[y=-4,x=-5, color=white, label=4]{E}
\Vertex[y=-4,x=-3, color=white, label=6]{F}
\Vertex[y=-4,x=-1, color=white, label=10]{G}
\Edge(A)(B)
\Edge(A)(C)
\Edge(B)(D)
\Edge(B)(E)
\Edge(C)(F)
\Edge(C)(G)

\Vertex[y=0,x=4, color=white, label=4]{A}
\Vertex[y=-2,x=2, color=white, label=1]{B}
\Vertex[y=-2,x=6, color=white, label=7]{C}
\Vertex[y=-4,x=5, color=white, label=6]{D}
\Vertex[y=-4,x=7, color=white, label=10]{E}
\Vertex[y=-4,x=1, color=white, label=1]{F}
\Vertex[y=-4,x=3, color=white, label=3]{G}
\Edge(A)(B)
\Edge(A)(C)
\Edge(C)(D)
\Edge(C)(E)
\Edge(B)(F)
\Edge(B)(G)
\end{tikzpicture}
}
\end{center}

The following proposition is easy to establish. Hence, the proof is left to the reader.
\Prp{20202020}{
Let $T(V,E)$ be a type-$2$ convex tree and $v$ be the root vertex of it. Then there exists at least one leaf node $v_{_{\ell}}$ such that the path $P$ connecting $v$ and $v_{_{\ell}}$ satisfies the weighted increasing monotonicity. In other words, if the path $P=\{v,v_1,\cdots,v_n,v_{_{\ell}}\}$, then the inequality
$w(v)\leq w(v_1)\leq\cdots\leq w(v_n)\leq w(v_{\ell})$
holds.
}

Before stating the results for the next section, we need to introduce a new terminology. The symbol $V_v$ will be used to denote the set consisting of all the descendant vertices of $v$, including the vertex $v$.\\

The next theorem demonstrates a stability result for an approximately type-$2$ convex tree. More precisely, it shows that if the vertices of the tree bear weights in an approximate monotonic arrangement, then the tree almost possesses type-$2$ convexity.\\
\Thm{11,000}{Let $\varepsilon>0$ and $T(V,E,w)$ be a rooted tree such that for any $v\in V$ with all its descendant vertices, it satisfies the following inequality
\Eq{11,000}{
w(v)\leq \inf_{v\in V_{_{v}}}w(v)+\varepsilon,
}
then there exists a type-$2$ convex tree $T(V,E,\widehat{w})$ such that $\|w-\widehat{w}\|_{\infty}\leq{\varepsilon}/2$ holds.
}
\begin{proof}
We will prove this theorem in three parts. First, we will show that for any tree $T(V,E,w'')$, it is possible to construct a type-$2$ convex minorant (possibly non-trivial). Let $v\in V$ be an internal vertex with atleast two child nodes and $v_1$ and $v_2$ are two of those. We consider the iterative functions $w_n:V\to\R_+$ as follows
\Eq{11002}{
w_0(v)=w''(v)\qquad \qquad \qquad 
w_n(v)=\min\bigg\{w_{_{n-1}}(v),\dfrac{w_{_{n-1}}(v_1)+w_{_{n-1}}(v_2)}{2}\bigg\}\quad (n\in\N).
}
From the construction, it is clear that $\Big(w_n(v)\Big)_{n\in\N}$ is a pointwise decreasing sequence bounded below by 0. This ensures its convergence. We assume, $\underset{n\to\infty}{\lim}w_n=\widehat{w}$ and thus for any internal vertex $v\in V$ together with its two child nodes $v_1$ and $v_2$ satisfies the following inequality
\Eq{*}{
\widehat{w}(v)=\lim_{n\to\infty} w_n(v)
&=\lim_{n\to\infty}\min\bigg\{w_{_{n-1}}(v),\dfrac{w_{_{n-1}}(v_1)+w_{_{n-1}}(v_2)}{2}\bigg\}\\
&\leq \lim_{n\to\infty}\dfrac{w_{_{n-1}}(v_1)+w_{_{n-1}}(v_2)}{2}\\
&=\dfrac{\widehat{w}(v_1)+\widehat{w}(v_2)}{2}.
}
This shows that the tree $T(V,E,\widehat{w})$ is a type-$2$ convex minorant of $T(V,E,w'')$.\\

Next, we are going to show that for any arbitrarily chosen internal vertex $v$, if the two weighted trees $T(V,E,w')$ and $T(V,E,w'')$ satisfy the following inequality
\Eq{11001}{
 w'(v)\leq\inf_{v\in V_{_{v}}}{w''(v)},
}
then there exists a type-$2$ convex tree $T(V,E,\widehat{w})$ such that the inequality
$w'(v)\leq \widehat{w}(v)\leq w''(v)$ holds for all $v\in V.$ To show the assertion, we consider $v\in V$ be arbitrary internal node and $v_1, v_2$ are two child nodes of it. Then, due to the convex combinations along with our assumption at \eq{11001}, the condition \eq{11002} can be extended to the following system of inequalities
\Eq{*}{
w'(v)&\leq \inf_{v\in V_{_{v}}}{w''(v)}\leq w_1(v)=\min\bigg\{w''(v),\dfrac{w''(v_1)+w''(v_2)}{2}\bigg\}\leq w''(v).\\
w'(v)&\leq \inf_{v\in V_{_{v}}}{w''(v)}\leq w_2(v)=\min\bigg\{w_1(v),\dfrac{w_1(v_1)+w_1(v_2)}{2}\bigg\}\leq w''(v).\\
\hdots \hdots& \hdots \hdots\hdots \hdots\hdots \hdots\hdots \hdots\hdots \hdots\hdots \hdots\hdots \hdots\hdots \hdots\hdots \hdots\hdots \hdots \hdots\\
\hdots \hdots& \hdots \hdots\hdots \hdots\hdots \hdots\hdots \hdots\hdots \hdots\hdots \hdots\hdots \hdots\hdots \hdots\hdots \hdots\hdots \hdots \hdots\\
w'(v)&\leq \inf_{v\in V_{_{v}}}{w''(v)}\leq w_n(v)=\min\bigg\{w_n(v),\dfrac{w_n(v_1)+w_1(v_2)}{2}\bigg\}\leq w''(v).\\
\hdots \hdots &\hdots \hdots\hdots \hdots\hdots \hdots\hdots \hdots\hdots \hdots\hdots \hdots\hdots \hdots\hdots \hdots\hdots \hdots\hdots \hdots \hdots\\
\mbox{so on}.
}
Since, $\Big(w_n(v)\Big)_{n\in\N}$ is a decreasing sequence bounded below by $\underset{{v\in V_{_{v}}}}\inf{w''(v)}$, this guarantees the convergence of the sequence. As shown in the first part of the proof, we will have a type-$2$ convex tree $T(V,E,\widehat{w})$ sandwiched in between the graphs $T(V,E,w')$ and $T(V,E,w'')$.\\

Finally, we can give the proof of the theorem as a corollary of the above result. We re-write the inequality \eq{11,000} as follows
\Eq{*}{
w(v)-\dfrac{\varepsilon}{2}\leq \inf_{v\in V_v} \bigg(w(v)+\dfrac{\varepsilon}{2}\bigg).
}
We consider $w':=\max\{w-\varepsilon/2\,,\,0\}$ and $w'':=w+\varepsilon/2$. Hence, the above inequality is equivalent to \eq{11001}. Thus, there exists a convex function $\widehat{w}:V\to\R_+$ such that the following inequality holds
\Eq{*}{
w(v)-\dfrac{\varepsilon}{2}\leq \widehat{w}(v)\leq w(v)+\dfrac{\varepsilon}{2}\quad\mbox{for all}\quad v\in V.
}
This establishes our result and completes the proof.
\end{proof}
This introductory paper on weighted convex graphs leaves several tempting questions and avenues for further research.
In the first section of this paper, the assumption of 
$n$-strict convexity on the graph $G(V,E,w)$ can be relaxed by imposing weaker conditions. For instance, most of the characteristic results remain valid even if we consider a weighted graph endowed with $n$-convexity, provided that no three adjacent vertices carry the same weight. Besides, in many results, the connectivity condition can be relaxed. Improved studies related to Ulam-type stability results can be carried out by assuming less strict conditions and implementing new innovative mathematical techniques. One may also investigate the possible integration of higher-order convexity with graphs. Moreover, it remains an open problem how to embed a convex sequence into an unweighted tree of the same cardinality such that the resulting structure turns into a weighted convex tree.\\

\section*{Statements and Declarations}

\noindent\textbf{Funding.} 
The first author received financial support from the \textit{``Ipar a Veszprémi Mérnökképzésért'' Foundation} for the research presented in this work.
\vspace{6pt}

\noindent\textbf{Use of AI.} 
Artificial intelligence tools were used solely for language polishing and grammatical correction; all scientific content and conclusions are entirely our own.

\vspace{6pt}
\noindent\textbf{Competing Interests.} 
The authors declare that there are no financial or non-financial competing interests relevant to the contents of this article.

\vspace{6pt}
\noindent\textbf{Ethics Approval.} 
Not applicable. This study does not involve human participants or animals.

\vspace{6pt}
\noindent\textbf{Consent to Participate.} 
Not applicable.

\vspace{6pt}
\noindent\textbf{Consent for Publication.} 
Not applicable.

\vspace{6pt}
\noindent\textbf{Data, Materials and/or Code Availability.} 
No datasets or code were generated or analysed during the current study.

\end{document}